\documentclass{paper}[12pt] 

\usepackage{amssymb}   
\usepackage{amsthm}
\usepackage{amsmath}
\usepackage{amsfonts}
\usepackage{latexsym} 
\usepackage{eucal}
\usepackage{indentfirst} 
\usepackage{graphicx} 
\usepackage{verbatim} 

\newtheorem{theorem}{Theorem}[section]

\newtheorem{lemma}[theorem]{Lemma}

\newtheorem{corr}[theorem]{Corollary}
\newtheorem{definition}[theorem]{Definition}

\def\R{{\mathbb R}} 
\def\Z{{\mathbb Z}} 
\def\N{{\mathbb N}} 

\begin{document}

\overfullrule=0pt
\baselineskip=24pt
\font\tfont= cmbx10 scaled \magstep3
\font\sfont= cmbx10 scaled \magstep2
\font\afont= cmcsc10 scaled \magstep2
\title{\tfont A Thermodynamic Classification of Real Numbers }
\bigskip
\bigskip
\author{Thomas Garrity\\  
Department of Mathematics and Statistics\\ Williams College\\ Williamstown, MA  01267\\ 
email:tgarrity@williams.edu }

\date{}

\maketitle
\begin{abstract}
 A new classification scheme for real numbers is given, motivated by ideas from statistical mechanics in general and work of Knauf \cite{Knauf1} and Fiala and Kleban \cite{Fiala-Kleban1} in particular. Critical for this classification of a real number will be the Diophantine properties of its continued fraction expansion.

\end{abstract}

\section{Introduction}

Though this paper is about number theory in general and about a 
classification scheme for real numbers in  particular, it has its roots in 
Thermodynamic Formalism, which was developed   in the 1960s by  Ruelle \cite{Ruelle2} \cite{Ruelle3} , Sinai \cite{Sinai1}  and others in an attempt to put statistical 
mechanics on a firm mathematical foundation.  Once done, the underlying 
mathematical scheme can then, in principle,  be applied to non-physical 
situations, using the original real-world interpretations to guide and 
influence what questions are to be asked and what structure is to be 
discovered.

This process has been begun in number theory.
In \cite{Knauf1}, Knauf developed a one-dimensional thermodynamic system based on the Farey fractions that exhibited phase transition.  In \cite{Fiala-Kleban1}, Fiala and Kleban generalized Knauf's work and showed  that their generalization has the same free energy as Knauf's.  We will  put these earlier works into a common linear algebra framework, allowing  us to make a seemingly minor, but actually significant, change in the original partition function.  We will produce, for each positive real number, a thermodynamic system.  Different real numbers will exhibit different free energies, giving us a new classification scheme for positive real numbers.  (This classification scheme can easily be extended to also include negative reals.)

In Section 2,  we give a brief overview of the parts of the statistical mechanics formalism that we will be using.  In particular, we will see the key importance of the partition function. In Section 3,  we tie this formalism to number theory, in particular to the Farey matrices.  In Section 3.3, we put Knauf's work into this language and do the same thing in Section 3.4 for Fiala and Kleban's work.  In Section 3.5, we show how to alter the earlier partition functions that will put us into the world of Diophantine analysis.  This is the section in which we are not just changing the notation from earlier work.  

In Section 4.1, we use our Diophantine partition function to give a new classification scheme for real numbers. In particular we develop the idea of a real number having a {\it  1-free energy limit}.  In Section 4.2, we show how this is naturally linked to continued fractions. The rest of Section 4 deals with proving that there are real numbers with 1-free energy limits, that there are reals without a 1-free energy limit, that all algebraic numbers have a $k$-free energy limits with $k>1$ and that all quadratic irrationals have 1-free energy limits.  We also show that $e$ has a $\sqrt{N}\log N$-free energy limit.  We will conclude with  open questions in Section 5.

There has been a lot of other work linking statistical mechanics to number theory.  There is other work of Knauf \cite{Knauf2} \cite{Knauf3} \cite{Knauf4}, of Guerra and Knauf \cite{Guerra-Knauf1}, of Contucci and Knauf \cite{Contucci-Knauf1}, of Fiala, Kleban and \"{O}zl\"{u}k \cite{Fiala-Kleban-Ozluk1}, of Kleban and \"{O}zl\"{u}k \cite{Kelban-Ozluk1}, of  Prellberg, Fiala and Kleban \cite{Prellberg-Fiala-Kleban1}, of Feigenbaum, Procaccia and Tel \cite{Feigenbaum-Procaccia-Tel1} and others.

 There is also the transfer operator method, applied primarily to the Gauss map, which allows, in a natural way, 
tools from functional  analysis to be used. We believe this was pioneered by Mayer  (see his \cite{Mayer2} for a survey) , and nontrivially extended by Prellberg \cite{Prellberg1}, by Prellberg and Slawny  \cite{Prellberg-Slawny1}, by Isola \cite{Isola1} and recently by Esposti, Isola and Knauf   \cite{Esposti-Isola-Knauf1}. An introduction to this work is in chapter nine of Hensley \cite{Hensley1}.  We will not be following this approach here.

As of August 2008, the web site 
\newline http://www.secamlocal.ex.ac.uk/people/staff/mrwatkin/zeta/physics.htm \newline
offers many other attempts over the years to find links between statistical mechanics and number theory.

Finally, I would like to thank Edward Burger for many interesting conversations about this work and Steven Miller and L. Pedersen for comments on an earlier draft. Also, I would like to thank  Peter Kleban,  Ali  \"{O}zl\"{u}k  and Thomas Prellberg for finding a significant error in an earlier draft.

\section{The Partition Function and the Free Energy}

This is a rapid fire overview of basic terms in statistical mechanics.  
For each $N \in \N$, we have a finite set $\mathcal{S}_N$, called the {\it state space}.  Let
$$E:\mathcal{S}_N\rightarrow \R^+$$
be a function that we call $\it{energy}$. The {\it partition function} is defined 
to be
$$Z_N(\beta)=\sum_{\sigma \in {\mathcal S}_N}e^{-\beta E(\sigma)}.$$
If we were modeling a physical system, the elements in the state space 
correspond to what can happen.  The variable $\beta$ corresponds to the 
inverse of the temperature. 
The underlying physical  assumption is that the probability that a 
system is in a state $\sigma \in \mathcal{S}_N$ will be
$$\mbox{Probability in state}\;\sigma = \frac{e^{-\beta 
E(\sigma)}}{Z_N(\beta)}.$$
While far from a proof, this interpretation makes sense, in that at  high 
temperatures (meaning for $\beta$ close to zero), all states become 
increasingly likely, while at low temperatures the most likely state 
increasingly becomes the state  with the lowest energy.

 There is a {\it free energy} if the following 
limit exists:

$$f(\beta)= \lim_{N\rightarrow \infty}\frac{\log(Z_N(\beta))}{N},$$
with the function $f(\beta)$ being called, naturally enough, the {\it free energy}.  It is 
believed that phase transitions occur at values of $\beta$ for which  $f(\beta)$ 
fails to be analytic. 

For almost all of this paper, our state space will be
$${\mathcal S}_N =\{ \sigma = (\sigma_1,\ldots,\sigma_N) : \sigma_i=0\;\mbox{or}\;1\}.$$
Thus each of our $\mathcal{S}_N$ will have order $2^N$.  We can think of 
our state space as having  $N$ site points, each having value $0$ 
or $1$.

The most famous example is the one-dimensional Ising model.  For 
convenience, we let each site have the value of $1$ or $-1$.  Thus for the 
Ising model, we have
$${\mathcal S}_N=\{\sigma=(\sigma_1,\ldots,\sigma_N):\sigma_i=\pm 1\}.$$
The energy function for the Ising model is
$$E(\sigma)=\sum_{i=1}^{N}\sigma_i \sigma_{i+1}.$$
Ising, in his 1925 thesis, showed that for this model there is no phase 
transition, meaning  he showed that  the free energy is an analytic function. For the
 two-dimensional analog, it is one of the great 
discoveries (originally by Onsager in 1944)  that phase transition does 
occur.    Most texts on statistical mechanics, such as  \cite{Thompson1}, describe 
the Ising model in detail.

Note that in the Ising model, a site will only interact with those other sites that 
are immediately adjacent to it.  This is an example of finite range 
interaction.  Since there is no phase transition for the one-dimensional Ising model, it 
was long believed that there would be no phase transition for any 
one-dimensional system.  But in the 1960s, it was discovered that phase 
transition can occur if the interactions are not of finite range but over 
possibly arbitrarily long distances.  A good introduction to this work is in Mayer's {\it The Ruelle-Araki transfer operator in classical statistical mechanics} \cite{Mayer1}.  Such interactions are called {\it long 
range interactions}.  In the following number theoretic models, it is key that the interactions are long range.

\section{Number Theoretic Partition Functions}

\subsection{General Set-up}
Fix a positive integer $k$.  For each positive integer $N$, our state space 
will be
$${\mathcal S}_N=\{(\sigma_1,\ldots, \sigma_N):\sigma_i=0,1,\ldots, k-1\}.$$
Thus ${\mathcal S}_N$ contains $k^N$ elements.

We define a new type of  product of an $N$-tuple of $n\times n$ matrices 
with an $M$-tuple of such matrices to be the $MN$-tuple:
$$(A_1,\ldots, A_N)(B_1,\ldots B_M)=(A_1B_1,A_1B_2, \ldots,A_NB_M).$$

For  matrices $A=(a_{ij})$ and $B=(b_{ij})$, we denote the Hilbert-Schmidt 
product (which is also called the Hadamard product) as
$$A*B=Tr(AB^T)=\sum_{1\leq i,j \leq n}a_{ij}b_{ij}.$$
For example, thinking of a $2\times 2$  matrix as an element of $\R^4$, then $A*B=Tr(AB^T)$ is simply the dot product of the two vectors.

Let $\mathcal{M}_n$ denote the space of $n\times n$ matrices. For a function
$$f:\mathcal{M}_n \rightarrow {\bf R}$$
and  for two $n\times n$ matrices $M$ and $A$, define
$$f(M)(\beta)|A=\frac{1}{|M*A|^{\beta}}f(MA^T),$$
following notation as in  \cite{Fiala-Kleban1}.
For $k$ $n\times n$ matrices 
$A_1,A_1,\ldots, A_{k}$, define
$$f(M)(\beta)|(A_1,\ldots, A_k) =\sum_{i=1}^k f(M)(\beta)|A_i.$$

Consider the map

$$Z:\N \times \mathcal{M}_n \times \R \times \mathcal{M}_n^k \times \Gamma( \mathcal{M}_n, \R) \rightarrow \R^*,$$
where $\N$ is the natural numbers, $\mathcal{M}_n$ is the space of $n\times n$ 
matrices, $\Gamma(\mathcal{M}_n, \R)$ is the space of functions from $n\times n$ matrices to the real 
numbers and $\R^*$ is the extended real numbers, defined by setting
$$Z(N,M,\beta,(A_0,\ldots,A_{k-1}),f)=f(M)(\beta)|(A_0,\ldots,A_{k-1})^N.$$
Here the notation $(A_0,\ldots,A_{k-1})^N$ is referring to the above newly defined product of tuples of matrices and hence can be viewed as short-hand for all products $A_{i_1}\cdots A_{i_N}$, with $0\leq i_l \leq k-1.$

We want to link this with partition functions. 
Fix an $n\times n$ matrix $M$ and also $k$ $n\times n$ matrices 
$A_0,A_1,\ldots, A_{k-1}$.
Let the function $f$ be the constant function $1$, or in words, let 
$f(B)=1(B)=1$
for all matrices $B\in \mathcal{M}_n$.
Then define the partition function to be
$$Z_N(M)(\beta)=Z(N,M,\beta,(A_0,\ldots,A_{k-1}),1)  =   1(M)(\beta)|(A_0,\ldots,A_{k-1})^N.$$

The ``physical" intuition is as follows.  Think of a one-dimensional lattice 
with $n$ sites.  At each site, there are $k$ possible states, each of 
which can either be indexed by a number $\sigma_i$ between $0$ and $k-1$ or 
by a matrix $A_{\sigma_i}$.  Then the states can be viewed as either all 
possible
$$\sigma=(\sigma_1,\ldots, \sigma_n)\in \Z/k\Z$$ 
or all matrices of the form
$$A_{\sigma_1}\cdots A_{\sigma_n}.$$
In order to get the above partition function, we set the energy  of a state
to be:
$$E(\sigma)=\log|M*(A_{\sigma_1} \cdots A_{\sigma_n})|.$$
In our applications, it  is more natural not to emphasize the energy 
function.

Using this language, we have natural recursion relations linking the 
partition function $Z_{n+1}$ with partition functions for various $Z_n$, 
with different choices for the matrix $M$.  More precisely

\begin{lemma}
$$Z_{N+1}(M)(\beta) = Z_N(MA_1^T)(\beta) +\ldots + Z_N(MA_k^T)(\beta).$$

\end{lemma}

This is a simple calculation.  It offers  a more general and natural form for  the 
key recursion 
relation (2) in \cite{Fiala-Kleban1}.  

\subsection{Farey Matrices}

This section continues the building of needed machinery.  (See, also,  section 4.5 in \cite{Graham-Knuth-Patashnik1}.  Another source would be \cite{Major1})

We develop  the Farey partitioning  of the extended real numbers.
Start with the set $${\cal F}_0=\left\{ \frac{1}{0},\frac{0}{1} \right\}$$
and define the Farey sum of two fractions, in lowest terms, to be
$$\frac{p}{q} \oplus \frac{r}{s} =\frac{p+r}{q+s}.$$
We now extend a given $n$th Farey set ${\cal F}_n$ to the $(n+1)$st Farey set 
by adding to ${\cal F}_n$ all of the 
terms obtained by applying the Farey sum.  Thus we have,
\begin{eqnarray*}
{\cal F}_0 &=& \left\{\frac{1}{0},\frac{0}{1}\right\} \\
{\cal F}_1 &=& \left\{\frac{1}{0},\frac{1}{1},\frac{0}{1}\right\} \\
{\cal F}_2 &=& \left\{\frac{1}{0},\frac{2}{1}, \frac{1}{1}, 
\frac{1}{2},\frac{0}{1}\right\} \\
{\cal F}_3 &=& \left\{\frac{1}{0},\frac{3}{1},\frac{2}{1}, \frac{3}{2},\frac{1}{1},
\frac{2}{3}, 
\frac{1}{2},\frac{1}{3},\frac{0}{1}\right\}.
\end{eqnarray*}

By reversing the order of the terms in ${\cal F}_n$, we get a partitioning 
of $\R^+ \cup \infty$.  Here we are thinking of $\frac{1}{0}$ as the 
point at infinity, which is why we are working with the extended real 
numbers $\R^*$.  

We now describe this partitioning in terms of iterations of matrix 
multiplication.  
Let $$A_0=\left( \begin{array}{cc} 1 &0 \\ 
1& 1 \end{array}\right) \; \mbox{and} \;\; A_1=\left( \begin{array}{cc} 1 &1 \\ 
0& 1 \end{array}\right).$$
These two matrices are key to understanding the Farey decomposition of the 
unit interval and  continued fraction expansions.  Further, these two 
matrices will  be key to the partition function of Knauf \cite{Knauf1}, of 
Fiala and Kleban \cite{Fiala-Kleban1}, and 
to our eventual use of partition functions for Diophantine approximations.

Note that 
$$\left( \begin{array}{cc} p &r \\ 
q& s \end{array}\right)A_0=\left( \begin{array}{cc} p+r &r \\ 
q+s& s \end{array}\right).$$
and
$$\left( \begin{array}{cc} p &r \\ 
q& s \end{array}\right)A_1=\left( \begin{array}{cc} p &p+r \\ 
q& q+s \end{array}\right).$$
Then ${\cal F}_1$, save for the point at infinity $\frac{1}{0}$,  can be obtained 
in the natural way by examining the  right columns of first $A_1$, then 
$A_0$.  Likewise, ${\cal F}_2$, save for the point at infinity $\frac{1}{0}$,  can be obtained 
in a similar  natural way by examining the the right columns of first $A_1^2$, then 
$A_1A_0$, then $A_0A_1$ and finally $A_0^2$.  In general the elements of
${\cal F}_n$, save of course for the  point $\frac{1}{0}$, are given by 
the
right columns of the $2^n$ products of matrices 
of $A_1$ and $A_0$.

These allow us to recover the continued fraction of a positive real number 
$\alpha$.  We know that  any  real number $\alpha$ can be written 

$$\alpha = a_0+\frac{1}{a_1+\frac{1}{a_2+\frac{1}{a_3 +\cdots }}},$$
which is usually denoted by 
$$\alpha = [a_0; a_1,a_2,a_3,\ldots],$$
where $a_0$ is an integer and the remaining $a_i$ are positive 
integers.  The 
number $\alpha$ is 
rational if and only if its continued fraction expansion terminates.  We 
say that the rational $\frac{p_m}{q_m}$ is the mth partial fraction for 
the number $\alpha$ if
$$\frac{p_m}{q_m}=[a_0;a_1,\ldots, a_m].$$

We want to use our Farey matrices to find the $n$th partial sum 
$\frac{p_n}{q_n}$ for a given number $\alpha$.  We return to our Farey 
numbers, but reverse the orders of the numbers:
\begin{eqnarray*}
{\cal F}_0 &=& \left\{\frac{0}{1},\frac{1}{0}\right\} \\
{\cal F}_1 &=& \left\{\frac{0}{1},\frac{1}{1},\frac{1}{0}\right\} \\
{\cal F}_2 &=& \left\{\frac{0}{1},\frac{1}{2}, \frac{1}{1}, 
\frac{2}{1},\frac{1}{0}\right\} \\
{\cal F}_3 &=& \left\{\frac{0}{1},\frac{1}{3},\frac{1}{2}, \frac{2}{3},\frac{1}{1},
\frac{3}{2}, 
\frac{2}{1},\frac{3}{1},\frac{1}{0}\right \}.
\end{eqnarray*}

We can thus view ${\cal F}_1$ as providing a splitting of the positive 
reals into two intervals: $[\frac{0}{1},\frac{1}{1}]$ and 
$[\frac{1}{1},\frac{1}{0}]$.  Note that  $[\frac{0}{1},\frac{1}{1}]$ can 
be thought as flipping the columns of $A_0$ and $[\frac{1}{1},\frac{1}{0}]$ 
can  be thought of as the flipping of the columns of $A_1$.

Likewise, ${\cal F}_2$ will split the positive reals into four intervals:
$[\frac{0}{1},\frac{1}{2}]$ (which can  be thought of as the flipping of the 
columns of $A_0A_0$), $[\frac{1}{2},\frac{1}{1}]$ (which can  be thought of as the flipping of the 
columns of $A_0A_1$), $[\frac{1}{1},\frac{2}{1}]$ (which can  be thought of as the flipping of the 
columns of $A_1A_0$) and $[\frac{2}{1},\frac{1}{0}]$ (which can  be thought of as the flipping of the 
columns of $A_1A_1$).  

This pattern continues. Consider the matrix 
 $$A_1^{a_0}A_0^{a_1}\cdots
A_1^{a_{N-1}}A_0^{a_N}$$  
with  $a_0\geq 0$ and $a_i>0$. Then the left column of the 
matrix will correspond to  
$[a_0;a_1,\ldots, a_N]$
while the right column will correspond to $[a_0;a_1,\ldots, a_{N-1}]$.

For matrices
 $$A_1^{a_0}A_0^{a_1}\cdots A_0^{a_{N-1}}
A_1^{a_N}$$
 with  $a_0\geq 0$ and $a_i>0$, then the left column will 
correspond to $[a_0;a_1,\ldots, a_{N-1}]$ while now the right column will $[a_0;a_1,\ldots, a_{N}]$.

Thus to determine the continued fraction expansion for a given positive real number $\alpha$, we just have to keep track of which interval  $\alpha$ is in for  a given $\mathcal{F}_N$.

\subsection{Knauf's work}
This section will show how Knauf's number theoretic partition function  \cite{Knauf1} 
can easily be put into the language of this paper.

Consider 
the  sets ${\cal F}_N \cap [0,1]$
\begin{eqnarray*}
{\cal F}_1 \cap [0,1] &=& \left\{\frac{0}{1},\frac{1}{1}  \right\} \\
{\cal F}_2  \cap [0,1]&=& \left\{\frac{0}{1},\frac{1}{2},\frac{1}{1}\right\} \\
{\cal F}_3  \cap [0,1]&=&\left \{\frac{0}{1},\frac{1}{3}, \frac{1}{2}, 
\frac{2}{3},\frac{1}{1}\right\} \\
{\cal F}_4 \cap [0,1] &=& \left\{\frac{0}{1},\frac{1}{4},\frac{1}{3}, \frac{2}{5},\frac{1}{2},
\frac{3}{5}, 
\frac{2}{3},\frac{3}{4},\frac{1}{1} \right\}.
\end{eqnarray*}
As is well known, as $N\rightarrow \infty$,  these sets will eventually 
contain all rational numbers in the unit interval $[0,1]$.  Then Knauf 
defines his partition function as 
$$Z_N^K(\beta) =\sum_{\frac{p}{q}\in {\cal F}_N\cap (0,1)} 
\frac{1}{q^{\beta}}.$$
(Note that  the $K$ is not being used as an  index but stands  for `Knauf'; the subscript $N$ is an index.)
In the limit we get 
$$Z^K(\beta)=\lim_{N\rightarrow \infty} Z_N^K(\beta) = 
\sum_{n=1}^{\infty}\frac{\phi(n)}{n^{\beta}},$$
where $\phi(n)$ is the Euler totient function.  
In turn,  $\sum_{n=1}^{\infty}\phi(n)/n^{\beta}$ is well-known to equal to 
$$\frac{\zeta(\beta -1)}{\zeta(\beta)},$$
for $\beta >2$, where
$\zeta(\beta)$ is the Riemann zeta function and is not defined for $\beta \leq 2$, showing that there is critical point phenomena 
for this one-dimensional system.  The  $\zeta(\beta -1)/\zeta(\beta)$ will show up throughout this paper and is why many of the later theorems are only true for $\beta >2$.

We now show how this can be put into the language of Section 4.1.  
Let 
$$M^K= \left( \begin{array}{cc} 0 &0 \\ 
0& 1\end{array}\right), $$
 and let $A_0$ and $A_1$ be the above Farey matrices.

Then the Knauf partition function is
$$Z_N^K(\beta)= 1(M^K)|A_0(A_0,A_1)^{N}.$$
The initial $A_0$ is just to insure that we are in the unit interval.  Also, for $A_0(A_0,A_1)^{N}$, we are using the new product for matrices defined in section 3.1 and not traditional matrix multiplication.

\subsection{Fiala-Kleban Work}

In \cite{Fiala-Kleban1}, Fiala and Kleban considered a different number 
theoretic partition function.  In the language of this paper, 
 let 
$$M^F= \left( \begin{array}{cc} 0 &0 \\ 
x& 1\end{array}\right).$$
(Again, the $F$ is not an index but stands for ``Fiala-Kleban".)
Then the new partition function will be
$$Z_N^F(x,\beta)= 1(M^F)|A_0(A_0,A_1)^{N}.$$
As mentioned earlier, the recursion relation (2) in Section 2 of \cite{Fiala-Kleban1}  is  
simply a special case of Lemma 3.1.

Two partition functions are said to have the same thermodynamics if their free energies are equal.
Thus,  the Knauf partition function $Z_N^K(\beta)$ and the Fiala-Kleban partition function
$Z_N^F(\beta)$ have the same thermodynamics if
$$\lim_{N\rightarrow \infty}\frac{\log (Z_N^K(\beta))}{N}=\lim_{N\rightarrow \infty}\frac{\log (Z_N^F(\beta))}{N}.$$
This equality is shown in Section 4 of  \cite{Fiala-Kleban1}, using as an intermediary tool a certain transfer operator and depending on the earlier work of Knauf \cite{Knauf1}.  
In \cite{Fiala-Kleban-Ozluk1},  Fiala,  Kleban and 
\"{O}zl\"{u}k showed the thermodynamics of the Knauf partition function is thermodynamically equivalent to 
a number of other number theoretic partition functions.  It is certainly the case, though, that thinking of the various matrices as vectors in ${\bf R}^4$ will yield more straightforward proofs of these equivalences.

\subsection{A Diophantine approach}
We now make a seemingly minor change in our choice for the matrix $M$ that 
will create a quite different thermodynamics, leading in the next section 
to a new classification of the real numbers. Set 

$$M= \left( \begin{array}{cc} 0 &-1\\ 
0& \alpha \end{array}\right)$$
for some real number $\alpha$.
Define the {\it Diophantine partition function} to be
$$Z_N(\alpha;\beta)= 1(M)|(A_0,A_1)^N.$$
Then we have
$$Z_N(\alpha;\beta) =\sum_{\frac{p}{q}\in {\cal F}_N}\frac{1}{|p-\alpha 
q|^{\beta}}.$$
Note that the terms that dominate the above sum occur when 
$|p-\alpha 
q|$
is small.  This places us firmly  in the realm of Diophantine analysis.  We will 
see that we can classify real numbers $\alpha$ by understanding the 
existence of free energy limits for the statistical systems associated to 
the partition function 
$Z_N(\alpha;\beta).$  Again, while this partition function is cast in the same overall 
language as Knauf and 
Fiala and Kleban, its thermodynamic properties will be quite different. This is what separates the present work from \cite{Knauf1} and \cite{Fiala-Kleban1}.

\section{Classifying real numbers via free energy}
\subsection{The classification}
With the notation above,
we have for each positive real number $\alpha$ and each positive integer $N$ the partition 
function
$Z_N(\alpha;\beta).$

\begin{definition} A real number $\alpha$ has a $k$-free energy limit if 
there is a number $\beta_c$ such that
$$\lim_{N\rightarrow \infty}\frac{\log(Z_N(\alpha;\beta))}{N^k}$$
exists for  all $\beta > \beta_c$
\end{definition}

For $k=1$, this is a number-theoretic version of the free energy of the system.
By an abuse of notation we will also say that $\alpha$ has a $f(N)$-free energy limit, for an increasing function $f(N)$ if
$$\lim_{N\rightarrow \infty}\frac{\log(Z_N(\alpha;\beta))}{f(N)}$$
exists for all $\beta >\beta_c$.

To see that this is a meaningful classification scheme for real numbers, we will establish the following three theorems:

\begin{theorem} There exists a  number that 
 has a 1-free 
 energy limit for  $\beta > 2$.
\end{theorem}

\begin{theorem} There exists a number that does not have a  1-free 
 energy limit, for any value of $\beta$.
\end{theorem}

 \begin{theorem}  Let $\alpha$ be a positive real number such that there is a positive constant $C$ and constant $d\geq 2$ with 
 $$\frac{C}{q^d}\leq |p-\alpha q|$$
 for all relatively prime integers $p$ and $q$.  Then $\alpha$
 has a k-free 
 energy limit for  $\beta > 2$, for any $k>1$, and in fact, the k-free energy limit is zero.
\end{theorem}
An easy consequence of the above is that all algebraic numbers  have  k-free 
 energy limits equal to zero, for any $k>2$. 

While we do not know if algebraic numbers have 1-free energy limits, we will show
\begin{theorem}
All quadratic irrationals have 1-free energy limits for $\beta >2.$
\end{theorem}

We will also show, in Section 4.7  that the number  $e$ has $\sqrt{N}\log(N)$-free energy limit.

\subsection{Links to continued fractions}
The goal of this paper is not just to give a new way to classify real numbers but  also to show how such a classification scheme follows from the thermodynamical formalism, fitting into a more  general framework.  But if all we wanted was the classification scheme, then it is possible to reframe our definitions so that there is no need for the language of statistical mechanics.  The goal of this section is to state the theorems that would allow us to avoid thermodynamics.  They will also be key to proving the five theorems of the previous section.

Let our positive real number $\alpha$ have continued fraction expansion \linebreak $[a_0;a_1,a_2,\ldots ]$.
We know that the {\it best rational approximations} to $\alpha$ are given by the rational numbers 
$[a_0;a_1,a_2,\ldots, a_m].$  The fractions $[a_0;a_1,a_2,\ldots, a_{m},k]$, with $0<k<a_{m+1}$, are called the {\it secondary convergents } to $\alpha$.  For a given $m$, we know that all the vectors corresponding to the  $[a_0;a_1,a_2,\ldots, a_{m-1},k]$, with $0\leq k<a_m$, are  on the same side of the line $x=\alpha y$, while $[a_0;a_1,a_2,\ldots, a_m]$  jumps to the other side of the line.

We now set some notation for the rest of the paper.  Given the positive real number $\alpha$,  for each positive integer $N$ there is associated a positive integer $m$ and nonnegative integer $k$ such that

$$N=a_0+a_1+\cdots + a_m +k,$$
with $0\leq k<a_{m+1}$.  We create  a subsequence of $\N$, denoted by $N_0,N_1,N_2,\ldots $ by setting
$$N_m=a_0+a_1+\cdots + a_m.$$
In this notation, given any $N\in \N$, we have
$$N_m\leq N < N_{m+1},$$
or, in other words,
$$a_0+a_1+\cdots + a_m \leq N=a_0+a_1+\cdots + a_m +k < a_0+a_1+\cdots + a_m +a_{m+1}.$$

Set 
$$\frac{p_N}{q_N}=[a_0;a_1,\ldots, a_{m},k],$$
with $0\leq k<a_{m+1},$
with  $p_N$ and $q_N$ having no nontrivial common factors.  We know that the fractions $\frac{p_{N_m}}{q_{N_m}}$ are the best rational approximations to the initial real $\alpha$. Denote
$$d_N=\frac{1}{| p_N-\alpha q_N |}.$$ We have the following chain of inequalities that will be key:

$$d_{N_{m-1}}<d_{N_{m}+1}<d_{N_{m}+2}<\cdots<d_{N_{m}+a_{m+1}-1}<d_{N_{m}},$$
which are well-known (for motiviation, see the chapter on continued fractions in \cite{Stark1}).

We will show 
\begin{theorem}
For any $\beta >2$ and for any positive real number $\alpha$, we have for all positive integers $N$ that
$$\frac{\log(d_N^{\beta})}{N}\leq \frac{\log(Z_N(\alpha;\beta))}{N}\leq \frac{\log 
(\frac{\zeta(\beta-1)}{\zeta(\beta)}Nd_{N_m}^{\beta})}{N}.$$
\end{theorem}

To show that there are numbers that do not have $1$-free energy limits, we will construct an $\alpha$ so that a subsequence of $\frac{\log(d_N^{\beta})}{N}$ approaches infinity.  In turn, to show that  that there are numbers that  have $1$-free energy limits, we will construct an $\alpha$ so that  $\frac{\log 
(\frac{\zeta(\beta-1)}{\zeta(\beta)}Nd_{N_m}^{\beta})}{N}$ approaches zero, forcing $\frac{\log(Z_N(\alpha;\beta))}{N}$ to approach zero (and thus guaranteeing that the limit exists).

\subsection{Preliminaries}

The reason why we look at the products of the matrices $A_0$ and $A_1$ is 
that the right columns will correspond to all of the integer lattice 
points 
$\bigl( \begin{smallmatrix}p \\ q \end{smallmatrix} \bigr)$ with $p$ and $q$ relatively prime in the the first quadrant of the plane.  If $A_{i_1}\cdots A_{i_n}$ has right column  $\bigl( \begin{smallmatrix}p \\ q \end{smallmatrix} \bigr)$, then
$$M*(A_{i_1}\cdots A_{i_n})=(-1,\alpha)\cdot  \begin{pmatrix}p \\ q \end{pmatrix}= \alpha q -p$$

Let
$$v_1=\begin{pmatrix}p_1 \\ q_1 \end{pmatrix}, v_2=\begin{pmatrix}p_2 \\ q_2 \end{pmatrix}$$
be vectors that satisfy
$\det (v_1,v_2) =\pm 1.$
Let $C(v_1,v_2)$ denote the cone of integer lattice points defined by:
$$C(v_1,v_2)=\{av_1+bv_2: a,b \; \mbox{relatively prime nonnegative integers}\}.$$
 Let $C_N(v_1,v_2)$ be the  subset of $C(v_1,v_2)$ consisting 
of vectors that are the right columns of all possible $A_{i_1}\cdots 
A_{i_N}$. Then set
$$Z_N(\alpha; \beta,v_1,v_2)=\sum_{\bigl( \begin{smallmatrix}p \\ q \end{smallmatrix} \bigr) \in C_N(v_1,v_2)} \frac{1}{|\alpha q- p|^{\beta}}.$$
In the same way, we set 
$$Z(\alpha; \beta,v_1,v_2)=\sum_{\bigl( \begin{smallmatrix}p \\ q \end{smallmatrix} \bigr)\in C(v_1,v_2)} \frac{1}{|\alpha q- p|^{\beta}}.$$

Suppose we have integer lattice vectors $v_1=\bigl( \begin{smallmatrix}0 \\ 1 \end{smallmatrix} \bigr), v_2,\ldots, v_m$ on one side of the line  $x=\alpha y$ and  integer lattice  vectors $w_1=\bigl( \begin{smallmatrix}1 \\ 0 \end{smallmatrix} \bigr), w_2,\ldots, w_n$ on the other  side of the line  $x=\alpha y$
such that $\det(v_i,v_{i+1})= -1$, $\det(w_i,w_{i+1})= 1$ and $\det(v_m,w_n)= 
-1$.  Then we have
$$Z_N(\alpha;\beta)\leq \sum_{i=1}^{m-1}Z_N(\alpha;\beta,v_i,v_{i+1}) +Z_N(\alpha;\beta,v_m,w_n) 
+\sum_{i=1}^{n-1}Z_N(\alpha;\beta,w_i,w_{i+1}).$$
We need the above to be an inequality since there is ``overcounting" on the right hand side, since, for example, the part of $Z_N(\alpha;\beta,v_{i-1},v_{i})$ coming from the vector $v_{i}$ also appears as a term in $Z_N(\alpha;\beta,v_i,v_{i+1})$. The key, as we will see, is that the  $Z_N(\alpha;\beta,v_m,w_n) $ term contributes the most to the partition 
function $Z_N(\alpha;\beta)$.   For the rest of this section, we want to  control the sizes of   the various 
$Z_N(\alpha;\beta,v_i,v_{i+1})$ and $Z_N(\alpha;\beta,w_i,w_{i+1})$.

We first return to the more general case of two vectors  $v_1=\bigl( \begin{smallmatrix}p_1 \\ q_1 \end{smallmatrix} \bigr), v_2=\bigl( \begin{smallmatrix}p_2 \\ q_2 \end{smallmatrix} \bigr)$  with $\det (v_1,v_2) = \pm 1$, under the additional assumption that $v_1$ and $v_2$ lie on the same side of the line $x=\alpha y$.
  
For 
$$d_1=\frac{1}{|p_1-\alpha q_1|} , d_2=\frac{1}{|p_2-\alpha q_2|},$$
suppose that $d_1<d_2$, which means that the line through the origin and the point $(p_2,q_2)$ is closer to the line $x=\alpha y$ than the line through the origin and $(p_1,q_1)$.

We want to show that 
$$Z(\alpha; \beta,v_1,v_2) <|d_2|^{\beta}\frac{\zeta(\beta -1)}{\zeta (\beta)}.$$

We know that 
$$\frac{1}{d_1}=\left|(-1,\alpha)\cdot \begin{pmatrix}p_1 \\ q_1 \end{pmatrix} \right|, \frac{1}{d_2}=\left|(-1,\alpha)\cdot  \begin{pmatrix}p_2 \\ q_2 \end{pmatrix} \right|.$$

Let $v$ be some 
integer lattice point in the cone $C(v_1,v_2)$.  Then there are relatively prime  positive 
integers $a$ and $b$ with $v=av_1+bv_2$.

We have 
\begin{eqnarray*}
(a+b) |(-1,\alpha)\cdot v_2| &<& | (-1,\alpha)\cdot (av_1+bv_2)|\\
&=& | (-1,\alpha)\cdot v|\\
&<& (a+b) |(-1,\alpha)\cdot v_1|.
\end{eqnarray*}

Inverting and raising everything to the power of $\beta$, we have

$$\frac{|d_1|^{\beta}}{(a+b)|^{\beta}}<\frac{1}{| (-1,\alpha)\cdot v|^{\beta}}<\frac{|d_2|^{\beta}}{(a+b)|^{\beta}}.$$

Summing over every vector in  $C(v_1,v_2)$, we get 
$$|d_1|^{\beta}\sum\frac{1}{|a+b|^{\beta}}<Z(\alpha;\beta,v_1,v_2) <|d_2|^{\beta}\sum\frac{1}{|a+b|^{\beta}},$$
where in the first and third summation we are summing over all relatively prime positive integers $a$ and $b$.  
It is well-known, as mentioned in Section 4.3 of Knauf \cite{Knauf1}, that,
$$\sum\frac{1}{|a+b|^{\beta}}= \frac{\zeta(\beta -1)}{\zeta 
(\beta)}.$$
It is here that the $\frac{\zeta(\beta -1)}{\zeta 
(\beta)}$ makes its critical appearance.
We have our desired inequality.

\subsection{Proof of Theorem 4.6}

First, the partition function $Z_N(\alpha;\beta)$ is the sum of many positive 
terms, including $d_N^{\beta}$. Thus we have the lower bound
$$\frac{\log(d_N^{\beta})}{N}\leq \frac{\log(Z_N(\alpha;\beta))}{N}.$$

Now for the upper bound.  We are at the Nth stage.  Letting $t\leq N$, expressing  $N$ as 
$N=a_0+\cdots +a_m +k$, with $0\leq k<a_{m+1}$ and using the notation from section 4.2, we assume that the vectors $\bigl(\begin{smallmatrix}p_{t-1}\\q_{t-1}\end{smallmatrix}\bigr)$ and $\bigl(\begin{smallmatrix}p_{t}\\q_{t}\end{smallmatrix}\bigr)$ lie on the same side of the line $(x=\alpha y)$. This is equivalent to there being some $s\leq m$ with
$a_0+\cdots +a_s \leq t-1 <t<a_0+\cdots +a_{s+1}$.  
We know that 
$$\det \begin{pmatrix}p_{t-1} &p_t \\ q_{t-1}& q_t \end{pmatrix} = \pm 1.$$
We have 
\begin{eqnarray*}
Z_N\left(\alpha;\beta;  \begin{pmatrix}p_{t-1}  \\ q_{t-1} \end{pmatrix} ,  \begin{pmatrix}p_t \\  q_t \end{pmatrix} \right)
& <&Z\left(\alpha;\beta; \begin{pmatrix} p_{t-1}  \\ q_{t-1} \end{pmatrix} ,  \begin{pmatrix}p_t \\  q_t \end{pmatrix}\right)\\
& <&\frac{\zeta(\beta -1)}{\zeta(\beta)} d_t^{\beta}\\
&<& \frac{\zeta(\beta 
-1)}{\zeta(\beta)} d_{N_m}^{\beta},
\end{eqnarray*}
using that $d_t\leq d_{N_m}.$

Since there will be $N$ such cones, we have

$$\frac{\log(Z_N(\alpha;\beta))}{N}\leq 
\frac{\log 
\left(N\frac{\zeta(\beta-1)}{\zeta(\beta)}d_{N_m}^{\beta}\right)}{N},$$
finishing the proof.

\subsection{Proof of Theorem 4.2}

We know that the best rational approximations to a real number $\alpha$ are the 
$$\frac{p_{N_m}}{q_{N_m}}=[a_0;a_1,\ldots, a_m].$$
It is well known that 
$$q_{N_{m+1}}=a_{m+1}q_{N_m}+q_{N_{m-1}}.$$
Further (as in Lemma 7.2 of \cite{Burger1})
$$a_{m+1}q_{N_m}\leq \frac{1}{|p_{N_m} -\alpha q_{N_m}|}=d_{N_m} \leq (a_{m+1}+2)q_{N_m}.$$

Our  goal in this section is to construct a real number $\alpha$ for which 
$$\lim_{N\rightarrow \infty}\frac{\log \left(N
\frac{\zeta(\beta-1)}{\zeta(\beta)}(d_{N_m})^{\beta}\right)}{N}=0,$$
for $\beta >2$, which by Theorem 4.5 will force 
 $\alpha$ to have a $1$-free energy limit. 
 
If

 $$ \lim_{N\rightarrow \infty}\frac{\log \left(N
\frac{\zeta(\beta-1)}{\zeta(\beta)}(d_{N_m})^{\beta}\right)}{N} $$
exists, then it 
equals 

$$ \lim_{N\rightarrow \infty}\frac{\log N
}{N} + \lim_{N\rightarrow \infty}\frac{\log 
\frac{\zeta(\beta-1)}{\zeta(\beta)}}{N} + \lim_{N\rightarrow \infty}\frac{\beta \log d_{N_m}}{N} .$$
 The first two terms  in the above certainly go to zero.  Thus we must construct a real $\alpha$ so that $\lim_{N\rightarrow \infty}\frac{\beta \log d_{N_m}}{N}=0.$  Recalling  our notation that
 $N_m=a_0+\cdots + a_m \leq N=N_m +k$,  with $0\leq k<a_{m+1},$  we have for each  $N$ 
 $$\frac{\log d_{N_m}}{N}\leq \frac{\log d_{N_m}}{N_m}.$$
 Thus all we have to do is construct an $\alpha$ so that $\lim_{N\rightarrow \infty}\frac{\beta \log d_{N_m}}{N_m}=0.$

 For $\alpha =[a_0;a_1,\ldots ]$, define the function $f(m)$ by setting
$$a_{m+1}=q_{N_m}^{f(m)}.$$
We have
\begin{eqnarray*}
\frac{\beta \log(d_{N_m})}{N} &\leq & \frac{\beta \log(d_{N_m})}{N_m}      \\
&\leq& \frac{\beta \log[(a_{m+1}+2)q_{N_m}] }{N_m}\\
&\leq& \frac{\beta\log(2a_{m+1}q_{N_m})}{N_m} \\
&=&\frac{\beta\log(2)}{N_m}+ \frac{\beta\log\left(q_{N_m}^{f(m)}q_{N_m}\right)}{N_m} \\
&=& \frac{\beta \log(2)}{N_m}+ \frac{\beta[f(m)+1]\log(q_{N_m})}{N_m}
\end{eqnarray*}
Since the first term in the last equation goes to zero as $N\rightarrow 
\infty$, we have, if the limits exist, that 
$$\lim_{N\rightarrow \infty}\frac{\beta \log( d_{N_m})}{N_m}
\leq \lim_{N\rightarrow \infty}\frac{\beta [f(m)+1]\log(q_{N_m})]}{N_m}$$

 We now start with a $q_0$ and a $q_1$ and a function $f(m)$ and use these to create our number $\alpha$,   Let $q_0=1$, $q_1=2$, and $f(m)=m$ for $m\geq1$.  Then define for $m\geq 2$,
$a_{m+1}=q_{N_m}^{f(m)}.$

Now
\begin{eqnarray*}
\log(q_{N_m})&=&\log( a_{m}q_{N_{m-1}}+q_{N_{m-2}}) \\
&\leq &\log(2a_mq_{N_{m-1}}) \\
&=& \log( 2q_{N_{m-1}}^{f(m-1)+1}) \\
&=& \log(2) + (f(m-1)+1)\log(q_{N_{m-1}})
\end{eqnarray*}
Then we have 

$$\frac{\beta[f(m)+1]\log(q_{N_m})}{N_m} \leq \frac{\beta[f(m)+1]\log(2)}{N_m}+ 
\frac{\beta[f(m)+1](f(m-1)+1)\log(q_{N_{m-1}})}{N_m}$$
We will now use that 
\begin{eqnarray*} 
N_m&=&a_0+\ldots + a_m\\
&>& a_m \\
&=&q_{N_{m-1}}^{f{(m-1)}}
\end{eqnarray*}
Then 
$$\frac{\beta[f(m)+1]\log(q_{N_m})}{N_m} < \frac{\beta(m+1)\log(2)}{q_{N_{m-1}}^{m-1}} + 
\frac{\beta(m+1)m\log(q_{N_{m-1}})}{q_{N_{m-1}}^{m-1}}$$
which has limit zero as $N \rightarrow \infty$.  Thus with the choice of 
the function $f(m)=m$ we have constructed  a real number that has 1-free energy limit, finishing the proof of Theorem 4.2.

 \subsection{Proof of Theorem 4.3}
 
 We use from Theorem 4.5 that 
 $\frac{\log(d_N^{\beta})}{N}\leq \frac{\log(Z_N(\alpha;\beta))}{N}.$  We will construct a real number $\alpha$ so that 
 $$\lim_{N_m\rightarrow \infty}\frac{\log(d_{N_m}^{\beta})}{N_m} = \infty.$$
 This will mean that $\lim_{N\rightarrow \infty} \frac{\log(Z_N(\beta))}{N}$ will not exist, which is the goal of Theorem 4.3.

 Proceeding as in the proof of Theorem 4.2, we define $\alpha =[a_0;a_1,a_2,\ldots]$ inductively on the $a_m'$s by 
 setting, as before,  $a_0=1$ and $a_{m+1}=q_{N_m}^{f(m)}$ but now defining  $f(m)$ as
$$f(m)=a_0+\cdots + a_m.$$
We use that $a_{m+1}q_{N_m} \leq d_{N_m}$ and $N_m=a_0+\cdots + a_{m}$, we have  
\begin{eqnarray*}
\frac{\log d_{N_m}^{\beta}}{N_m}&\geq &\frac{\log \left(a_{m}q_{N_m}\right)^{\beta}}{N_m}\\
&=&\frac{\beta \log \left(q_{N_m}^{f(m)}q_{N_m}\right)}{N_m} \\
&=&\frac{\beta(f(m)+1)\log q_{N_m}}{N_m}\\
&=&\frac{\beta(a_0+\cdots a_{m}+1)\log q_{N_m}}{a_0+\cdots + a_{m}}.
\end{eqnarray*}
Since the last term goes to infinity as $N_m\rightarrow \infty$, we are done.

\subsection{Proof of Theorem 4.4}

We assume that $\alpha$ is a positive real number such that there is a positive constant $C$ and a constant $d\geq 2$ with 
 $$\frac{C}{q^d}\leq |p-\alpha q|$$
 for all relatively prime integers $p$ and $q$. 

In particular, we have 
$$\frac{C}{q_{N_m}^d}\leq |\alpha-\frac{p_{N_m}}{q_{N_m}}|.$$
In our
notation, this is 
$$d_{N_m}\leq cq_{N_m}^{d-1},$$
where $c$ is a constant depending on $\alpha$ but not on ${N_m}$. 
In a similar argument as in  Theorem 4.2, the number $\alpha$ will have a k-free energy limit 
if we can show 
that 
$$\lim_{N_m \rightarrow \infty}\frac{\log(d_{N_m})}{{(N_m)^k}}=0.$$
This will certainly happen if we can show
$$\lim_{N_m \rightarrow \infty}\frac{\log(cq_{N_m}^{d-1})}{{(N_m)^k}} =0$$
and hence if
$$\lim_{N_m\rightarrow \infty}\frac{\log(q_{N_m})}{{(N_m)^k}}=0.$$

We know that
$$q_{N_m}= a_{m}q_{N_{m-1}} + q_{N_{m-2}}.$$
 Then we get 
 \begin{eqnarray*}
 q_{N_m} &=& a_{m}q_{N_{m-1}} + q_{N_{m-2}}\\
&\leq & 2a_{{m}}q_{N_{m-1}} \\
&=& 2a_{{m}}(a_{m-1}q_{N_{m-2}} + q_{N_{m-3}}) \\
&\leq & 2^2a_{m}a_{m-1}q_{N_{m-2}}\\
&\vdots & \\
&\leq & 2^{m+1}a_{m}a_{m-1}\dots a_0 .
\end{eqnarray*}

Thus we want to examine
$$\lim_{N_m\rightarrow \infty}\frac{\log(2^{m+1}a_{m}a_{m-1}\dots a_0  )}{(N_m)^k}$$
or

$$\lim_{N_m\rightarrow \infty}\frac{(m+1)\log(2) +  \log(a_0) + \ldots +
\log(a_{m})}{(a_0+\ldots + a_{m})^k}.$$
Since it is always the case that  $a_0+\ldots + a_{m} \geq m+1,$ the above limit must always be
 zero for any $k>1$. 
 
 \begin{corr}  All positive algebraic numbers of degree greater than or equal to two have k-free energy limits equal to zero, for any $k>1$.
 \end{corr}
 
Liouville's Theorem (see Theorem 1.1 in \cite{Burger-Tubbs1}) states that all irrational algebraic numbers $\alpha$ have the property there is a constant $d\geq 2$ with 
 $$\frac{C}{q^d}\leq |p-\alpha q|$$
 for all relatively prime integers $p$ and $q$, allowing us to immediately use the above theorem.

  \subsection{Quadratic Irrationals have 1-free energy limits}
  The key will be that every quadratic irrational has an eventually periodic continued fraction expansion (see section 7.6   in \cite{Stark1}).  We will show first, for a quadratic irrational, that  $\lim_{m\rightarrow \infty} \log q_{N_m}/N_m$ exists and then show that the existence of this limit is equivalent to the number   having a 1-free energy limit.

  Let $\alpha$ be a quadratic irrational.    We can write
$$\alpha =[b_0:,b_1,\ldots , b_p,c_1, c_2, \ldots, c_l, c_1, c_2, \ldots, c_l, \ldots].$$
The period length is $l$.
For notational convenience, set
\begin{eqnarray*}
b&=&b_0 +\cdots + b_p \\
c&=&c_1+\cdots + c_l \\
B&=&\left(\begin{array}{cc} 0&1\\  1&b_p \end{array}\right)\cdots \left(\begin{array}{cc} 0&1\\  1&b_0 \end{array}\right)\\
C&=&\left(\begin{array}{cc} 0&1\\  1&c_l \end{array}\right)\cdots \left(\begin{array}{cc} 0&1\\  1&c_1 \end{array}\right).
\end{eqnarray*}
For each positive $m>p$, there are unique integers $k\geq 0$ and $n$  with $0\leq n<l$ such that
$$m=p+kl+n.$$
For this $n$, set
\begin{eqnarray*}
d_n&=&c_1+\cdots + c_n\\
D_n&=&\left(\begin{array}{cc} 0&1\\  1&c_n \end{array}\right)\cdots \left(\begin{array}{cc} 0&1\\  1&c_1 \end{array}\right).
\end{eqnarray*}
Note that there are only a finite number of possibilities for the $d_n$ and the $D_n$.

From section 7.6 of Stark, we have
$$\left(\begin{array}{cc} q_{N_{m-1}}&p_{N_{m-1}}\\  q_{N_m}&p_{N_m} \end{array}\right)=D_nC^kB.$$
 Now 
$N_m=b+kc +d_n$
 grows linearly with respect to $k$.  Since $  q_{N_m}$ is the lower left term of the matrix $D_nA^kB$, we also have that $\log  q_{N_m}$ grows linearly with respect to $k$, with constant leading coefficient.  Thus $\lim_{m\rightarrow \infty}\frac{\log q_{N_m}}{N_m}$ must exist.

  Now to show why this will imply that $\alpha$ has a 1-free energy limit.
  Since
   $d_{N_{m-1}}\leq d_N\leq d_{N_m}$ and since $N_m\leq N \leq N_{M+1}$, we have from Theorem 4.6 that
    $$\frac{\log(d_{N_{m-1}}^{\beta})}{N_{m+1}}\leq \frac{\log(Z_N(\alpha;\beta))}{N}\leq \frac{\log 
(\frac{\zeta(\beta-1)}{\zeta(\beta)}Nd_{N_m}^{\beta})}{N_m} .$$
From our earlier work, we know that $\alpha$ will have a 1-free energy limit, for $\beta >2$, if we can show

$$\lim_{m\rightarrow \infty}\frac{\log d_{N_{m-1}}}{N_{m+1}}=\lim_{m\rightarrow \infty}\frac{\log d_{N_{m}}}{N_{m}}.$$
We will see that these limits will exist if and only if $\lim_{m\rightarrow \infty}\frac{\log q_{N_m}}{N_m}$ exists.

Let $\alpha =[a_0; a_1, a_2, \ldots ]$.  We know that 
$N_{m+1}=N_{m-1}+a_m + a_{m+1}. $  For $m>p$, $a_m$ and $a_{m+1}$ can only be chosen from  finite  set $\{c_1,\ldots c_l\}$.  Thus 
$$\lim_{m\rightarrow \infty}  \frac{\log d_{N_{m-1}}}{N_{m+1}} =\lim_{m \rightarrow \infty}   \frac{\log d_{N_{m-1}}}{N_{m-1}}$$
since the denominators, for all $m$, differ by a fixed amount (and of course since $N_m\rightarrow \infty$).

We know that 
$a_{m+1}q_{N_m}\leq d_{N_m}\leq (a_{m+1} +2)q_{N_m}$.
 Then, if the limits exist,  we have
$$\lim_{m \rightarrow \infty}  \frac{\log a_{m+1}q_{N_m}}{N_m}\leq \lim_{m \rightarrow \infty}  \frac{\log d_{N_m}}{N_m}\leq \lim_{m \rightarrow \infty}  \frac{\log (a_{m+1}+2)q_{N_m}}{N_m}.$$
Since the $a_{m+1}$ are bounded and the $N_m\rightarrow \infty$, we have
$$\lim_{m\rightarrow \infty} \frac{\log a_{m+1}}{N_m} =  \lim_{m\rightarrow \infty} \frac{\log (a_{m+1}+2)}{N_m} =0,$$ giving us that $$\lim_{m \rightarrow \infty}  \frac{\log d_{N_m}}{N_m}= \lim_{m \rightarrow \infty}    \frac{\log q_{N_m}}{N_m},$$
concluding the proof.

\begin{corr}
The golden ratio
$\phi=\frac{1+\sqrt{5}}{2}$ has 1-free energy limit equal to $\log \phi,$ for $\beta>2$.
\end{corr}  

For background on $\phi$, see Chapter 1.7 in \cite{Rockett-Szusz1}.  The key is that continued fraction expansion for $\phi$ is
$$\phi=[1:1,1,1,1,\ldots ].$$
From the above theorem, we know that the 1-free energy limit, for $\beta >2$ is
$$\lim_{m\rightarrow \infty}\frac{\log q_{N_m}}{N_m}.$$
For $\phi$, we know that 
$N_m=m+1$ and that $q_{N_m} = F_m$, where $F_m$ is the $m$th Fibonacci number. We know that 
$$F_m=\frac{1}{\sqrt{5}}\left(\left(\frac{1+\sqrt{5}}{2}\right)^m - \left(\frac{1-\sqrt{5}}{2}  \right)^k      \right).$$
Then the  1-free energy limit will be
\begin{eqnarray*}
\lim_{m\rightarrow \infty}\frac{\log q_{N_m}}{N_m}&=& \lim_{m\rightarrow \infty}\frac{\log F_m}{m+1} \\
&=&\log \phi,
\end{eqnarray*}
as desired.

Of course, similar arguments can be used to find the 1-free energy limits for many quadratic irrationals.

\subsection{The number $e$ has $\sqrt{N}\log(N)$-free energy limit}
We will actually show this for $e-1$.  The continued fraction expansion for $e-1$ is
$$[1;1,2,1,1,4,1,1,6,\ldots].$$
Thus for $e-1$, we have
$$a_{3k}=a_{3k+1}=1\; \;\mbox{and} \; \;a_{3k+2}=2(k+1).$$ Let  $q_{N_m}$ denote the denominators of the fraction $[a_0;a_1,\ldots a_m]$ associated to $e-1$, where $N_m=a_0 + \cdots + a_m.$

In an analogous fashion to the proof that quadratic irrationals have $1$-free energy limits, we will first show for $e-1$ that
$$\lim_{N_m\rightarrow \infty}\frac{\log q_{N_m}}{\sqrt{N_m}\log(N_m)}$$
exists and then show that the existence of this limit will give us our result.

 We know that 
 $$a_{m}q_{N_{m-1}}<q_{N_m}= a_{m}q_{N_{m-1}} + q_{N_{m-2}}.$$
 Then we have 
 $$ a_0a_1\cdots a_m  <  q_{N_m}\leq  2^{m+1}a_{m}a_{m-1}\dots a_0.$$
Then, provided the limits exist, we have

 $$\lim_{N_m\rightarrow \infty}\frac{\log  a_0a_1\cdots a_m}{\sqrt{N_m}\log(N_m)}\leq \lim_{N_m\rightarrow \infty}\frac{\log q_{N_m}}{\sqrt{N_m}\log(N_m)}\leq \lim_{N_m\rightarrow \infty}\frac{\log 2^{m+1}a_{m}a_{m-1}\dots a_0 }{\sqrt{N_m}\log(N_m)}.$$
 
 We will show that 
 $$\lim_{N_m\rightarrow \infty}\frac{\log  a_0a_1\cdots a_m}{\sqrt{N_m}\log(N_m)}= \lim_{N_m\rightarrow \infty}\frac{\log 2^{m+1}a_{m}a_{m-1}\dots a_0 }{\sqrt{N_m}\log(N_m)},$$
 which will force $\lim_{N_m\rightarrow \infty}\log q_{N_m}/\sqrt{N_m}\log(N_m)$ to exist.
 
 We now have to look at the explicit values for the $a_m$.
There are three cases, depending on if $m$ is $0,1$ or $2$ mod 3.  We will let
$m=3k+2$ and show that the above limits exist and are equal.  The other two cases are similiar.
The above denominator is
$$1+1+2+1+1+4+1+\ldots + 2(k+1)= 2(k+1) + (k+1)(k+2),$$
which has quadratic growth in $k$.  The numerator on the left hand side of the above limit is
\begin{eqnarray*} 
\log a_0 +\cdots \log a_m&=& \log 1 + \log 1 +\log 2 + \log 1 + \log 1 +\log 4 +\cdots \log 2(k+1) \\
&=& k\log 2 + \log 1 + \log 2 +\cdots \log (k+1),
\end{eqnarray*}
 which has $k\log k$ growth. 
 Similarly for the right hand side of the above limit, we have
 $$\log 2^{m+1}a_{m}a_{m-1}\dots a_0 =\log(a_{m+1}+2) +(m+1)\log(2) +  \log(a_0) + \ldots +
\log(a_{m}),$$
which also has $k\log k$ growth, with the same leading coefficient.  The above limits will exist. 

Now to show that 
$$\lim_{N\rightarrow \infty} \frac{\log(Z_N(\alpha;\beta))}{\sqrt{N}\log N}$$
exists.
 
 In an argument similar to the proof of  Theorem 4.6, we have 
    $$\frac{\log(d_{N_{m-1}}^{\beta})}{\sqrt{N_{m+1}}\log N_{m+1}}\leq \frac{\log(Z_N(\alpha;\beta))}{\sqrt{N}\log N}\leq \frac{\log 
(\frac{\zeta(\beta-1)}{\zeta(\beta)}Nd_{N_m}^{\beta})}{\sqrt{N_m}\log N_m} .$$
Following the proof in the last section, we see that 
$$\lim_{m\rightarrow \infty}\frac{\log d_{N_{m-1}}}{\sqrt{N_{m+1}}\log N_{m+1}}=\lim_{m\rightarrow \infty}\frac{\log d_{N_{m-1}}}{\sqrt{N_{m-1}}\log N_{m-1}},$$
since $N_{m+1}-N_{m-1}$ grows at most linearly with respect to $m$ and since each $N_m$ grows quadratically with respect to $m$.

Since
$a_{m+1}q_{N_m}\leq d_{N_m}\leq (a_{m+1} +2)q_{N_m}q_{N_m}$, if the limits exist,  we have
$$\lim_{m \rightarrow \infty}  \frac{\log a_{m+1}q_{N_m}}{\sqrt{N_m}\log N_m}\leq \lim_{m \rightarrow \infty}  \frac{\log d_{N_m}}{\sqrt{N_m}\log N_m}\leq \lim_{m \rightarrow \infty}  \frac{\log (a_{m+1}+2)q_{N_m}}{\sqrt{N_m}\log N_m}.$$
Since the $a_{m+1}$   grow at most linearly with respect to $m$ and the $N_m$ grow quadratically, we have
$$\lim_{m\rightarrow \infty} \frac{\log a_{m+1}}{\sqrt{N_m}\log N_m} =  \lim_{m\rightarrow \infty} \frac{\log (a_{m+1}+2)}{\sqrt{N_m}\log N_m} =0,$$ giving us that $$\lim_{m \rightarrow \infty}  \frac{\log d_{N_m}}{\sqrt{N_m}\log N_m}= \lim_{m \rightarrow \infty}    \frac{\log q_{N_m}}{\sqrt{N_m}\log N_m},$$
concluding the proof.

\section{Conclusion}  

There are a number of directions for future work.  It would certainly be interesting to know if there are algebraic numbers besides the quadratics that have 1-free energy limits.  Also, what type of free energies are there for various other real numbers.  We have shown that  some real numbers have 1-free energy limits while others do not.  This too is just a beginning.  What types of limits are possible?  Are there numbers $\alpha$ for which the sequence $\frac{\log\left(Z_N(\alpha;\beta)\right)}{N}$ has any possible limit behavior?  For example, certainly we should be able to rig $\alpha$ so that any number can be the limit of the sequence.  In fact, we should be able to find such sequences with two accumulation points, three accumulation points, etc.  All of these should provide information about the initial real number $\alpha$.

Once we know that there is a free energy limit, then the most pressing question is to find for which values of $\beta$ is the free-energy non-analytic.  These points will be the analogs to critical point phenomena in physics.  For the Knauf approach, there have been a number of good papers (such as in \cite{Knauf1} \cite{Knauf5}  \cite{Knauf2} \cite{Knauf3} \cite{Knauf4} \cite{Boca1}  \cite{Contucci-Knauf1}  \cite{Esposti-Isola-Knauf1} \cite{Fiala-Kleban-Ozluk1} \cite{Guerra-Knauf1} \cite{Isola1} \cite{Kallies-Ozluk-Peter-Syder1} \cite{Manfred1}       ) exploring the nature of these critical points. For our Diophantine partition functions, these types of questions seem to be equally subtle.

The various partition functions $Z_N(\beta)$ depend on the choice of the two-by-two matrix $M$.  Different choices of $M$ give rise to different thermodynamics.   Knauf's $M^K$ and Fiala's and Kleban's  $M^F$  are of the form
$$\left( \begin{array}{cc} 0 &0\\ 
x& 1 \end{array}\right)$$
while the $M$ introduced in this paper is
$$\left( \begin{array}{cc} 0 &-1\\ 
0& \alpha \end{array}\right).$$
But certainly other $M$ can be chosen and studied.

Also, $Z_N(\beta)$ depends on the choice of the matrices $A_0$ and $A_1$.  Other choices for these matrices will lead to different thermodynamical systems, each with its own number-theoretic implications. 

Of course, why stick to two-by-two matrices.  This leads to multi-dimensional continued fractions.  There are many different multi-dimensional continued fraction algorithms.  See Schweiger \cite{Schweiger1} for a sampling of some.  Also, Major \cite{Major1} has some preliminary work on this.   Each of these will give rise to a thermodynamical system, again with meaning in number theory.

 Finally, there is the question of putting these results into the language of transfer operators. (See \cite{Ruelle1}, \cite{Mayer1}, \cite{Baladi1} for general references.)

\end{document}